\documentclass[12pt]{amsart} 
\usepackage{amsmath,amssymb,amscd}
\input xy
\xyoption{all}

\DeclareMathOperator{\rad}{rad}

\newtheorem{thm}{Theorem}[subsection]
\newtheorem{cor}[thm]{Corollary}
\newtheorem{lem}[thm]{Lemma}
\newtheorem{prop}[thm]{Proposition}
\theoremstyle{definition}
\newtheorem{defin}[thm]{Definition}
\newtheorem{ex}[thm]{Example}
\newtheorem{con}[thm]{Construction}
\newtheorem{alg}[thm]{Algorithm}
\newtheorem{rmk}[thm]{Remark}

\newcommand{\AAA}{\mathcal{A}}
\newcommand{\HHH}{\mathcal{H}}
\newcommand{\IA}{I_{\mathcal{A}}}
\newcommand{\LH}{L_{\mathcal{H}}}
\newcommand{\AAH}{A_{\mathcal{A},\,\mathcal{H}}}
\newcommand{\AaAH}{A^*_{\mathcal{A},\,\mathcal{H}}}

\newcommand{\BAH}{B_{\mathcal{A},\,\mathcal{H}}}

\newcommand{\BAxH}{B_{\mathcal{A}}}

\begin{document}
\author[Bj\"orner]{Anders Bj\"orner}
\address{Department of Mathematics, Royal Institute of Technology, S-100~44
Stockholm, Sweden}
\email{bjorner@math.kth.se}

\author[Peeva]{Irena Peeva}
\address{Department of Mathematics, Cornell University, Ithaca, NY 14853, USA}
\email{irena@math.cornell.edu}

\author[Sidman]{Jessica Sidman}
\address{Department of Mathematics and Statistics, 415 A Clapp Lab, Mount
Holyoke College, South Hadley, MA 01075, USA}
\email{jsidman@mtholyoke.edu}

\title[Arrangements defined by products of linear forms]
{Subspace arrangements defined by products of linear forms}


\begin{abstract}
We consider the vanishing ideal of an  arrangement of
linear subspaces in a vector space and investigate when this ideal can be
generated by products of linear forms.  We introduce a  combinatorial
construction (blocker duality) which yields  such generators 
in cases with a lot of combinatorial structure,
and we present  the examples that motivated our work.
We give a construction which
produces all elements of this type in the vanishing
  ideal of the arrangement. This leads to  an algorithm for deciding if
  the ideal is generated by products of linear forms.  We also
consider generic arrangements of points in ${\bf P}^2$ and lines
in ${\bf P}^3.$
\end{abstract}
\maketitle

\section{Introduction}

\noindent Throughout the paper $k$ is an infinite field. 
We consider an arrangement
$\mathcal{A}$ of $r$ linear subspaces in $k^n$; we assume that none of the
subspaces contains another.
Let $I_1,   \dots ,I_r$ be the linear ideals in $k[x_1,\dots ,x_n]$
that are the defining
ideals of the subspaces in $\mathcal{A}$.
Denote by $V_{\mathcal{A}}$  the union of the subspaces in $\mathcal{A}$. The
vanishing ideal of $V_{\mathcal{A}}$ is the reduced ideal $$I_{\mathcal{A}}=
I_1\cap  \dots \cap I_r\, .$$
The ideal defining a subspace arrangement arises in connection with topics as
diverse as independence numbers of graphs and
graph coloring (see \cite{li-li1}, \cite{li-li2}, \cite{lovasz},
\cite{deloera}, \cite{domokos}),
invariant theory \cite{derksen}, and symmetric function
theory \cite{haiman}.

  When $\mathcal{A}$ is an arrangement of
hyperplanes its
vanishing ideal $I_{\mathcal{A}}$ is a very simple object --  a principal
  ideal
generated by the product of linear forms that define the hyperplanes.
In general, the ideal $I_{\mathcal{A}}$ is generated by  products of
linear forms   up to a radical, since
  $\rad (I_1 \cdots
I_r)= \rad (I_1) \cap \dots \cap \rad (I_r)=
I_1 \cap \dots \cap I_r=
I_{\mathcal{A}} $,
but it is difficult to construct a nice system of
generators of $I_{\mathcal{A}}$ itself.
Geometrically, finding generators of $I_{\mathcal{A}}$
is related to detecting
low-degree hypersurfaces intersecting in $V_{\mathcal{A}}$.
   We will present examples where the ideal
$I_{\mathcal{A}}$ is generated by products of
linear forms
in many cases in which $\mathcal{A}$ has a great deal of combinatorial
structure.

We say that an ideal is {\it pl-generated} if it is generated by
products of linear forms.
In this paper, we study combinatorial properties of $\AAA$ that
are related to $I_{\mathcal{A}}$ being pl-generated.
We present the combinatorial point of view in \S 3 and the
ideal-theoretic point of view in \S 4.
The last section is entirely different in flavor:
in \S 5 we  study when $I_{\mathcal{A}}$ is pl-generated for arrangements in
${\bf P}^2$
and ${\bf P}^3$.
\smallskip

  In  \S 3 we introduce the notion of blocker duality, a
combinatorial  operation which, given a subspace arrangement $\AAA$
and an embedding $\mathcal{A} \subseteq \mathcal{H}$ into a hyperplane
arrangement $\HHH$,  produces another ``dual''  subspace arrangement $\AAA^*$.
It is not in general true that $\AAA^{**}=\AAA$, only that
$V_{\AAA^{**}} \supseteq V_{\AAA}$.

We provide an overview of the examples that motivated
this construction in \S \ref{sec: motivation}.  In \S
\ref{sec: blockers} we define blocker duality and demonstrate its basic
properties.  In \S \ref{sec: blocker ideal} we use blocker
duality to define a
pl-generated ideal $\BAH$
which is contained in $I_{\mathcal{A}}$.
Over an algebraically closed field we show that
$\AAA=\AAA^{**}$ if and only if $\rad (\BAH)
= I_{\mathcal{A}}.$

The stronger property, that $\BAH=
I_{\mathcal{A}}$, holds for our motivating examples, as well as for some other
fundamental examples which we discuss in \S \ref{sec: examples}.
One would like to determine a
combinatorial property of $\AAA$, viewed as an antichain
in the intersection lattice $L_{\mathcal{H}}$ of
$\mathcal{H}$, that makes it possible  to detect if $\BAH$ is radical
and explains the examples.
\smallskip

In \S 4, we depart from the beautiful examples where blocker duality
works and we
consider the more general situation where the blocker ideal
$B_{\AAA,\HHH}$ may fail to be
equal to
$I_{\AAA}$. This could be caused by the following two problems:
\begin{enumerate}
\item[(1)] $B_{\AAA,\HHH}$ may fail to capture all products of
linear forms in $\mathcal{H}$ that are contained in $I_{\mathcal{A}}.$
\item[(2)] It might not be possible to generate enough products of linear forms 
using only linear forms from $\mathcal{H}$.
\end{enumerate}

In \S 4.1 we solve the first problem by introducing the ideal
$F_{\mathcal{A},\mathcal{H}}$ which is larger than $B_{\AAA,\HHH}$.
It is constructed combinatorially, but  it
is also a natural algebraic object: $F_{\mathcal{A},\mathcal{H}}$ is 
the largest
ideal inside
$I_{\mathcal{A}}$ that is generated by products of linear forms in
$\mathcal{H}$.

In \S 4.2 we solve the second problem.
We prove that any given  embedding ${\mathcal{A}}
\subseteq
{\mathcal{H}}$ can be enlarged to an embedding ${\mathcal{A}}\subseteq
\tilde{\mathcal{H}}$ so that $F_{\AAA,\tilde\HHH}$ is the  ideal
    generated by all  products of linear forms inside $I_{\mathcal{A}}$. In
particular, Theorem \ref{prop: full ideal} shows that
if $\AAA$ has the pl-property then a system of  generators of
$I_{\mathcal{A}}$ that are  products of
linear forms can be constructed by a combinatorial procedure. As an  immediate
consequence we obtain  Algorithm \ref{alg}, which makes it possible to check
by computer whether a given  ideal
$I_{{\mathcal{A}}}$ is pl-generated.  For a related result see Proposition 1.1 
in \cite{li-li2}. 
\smallskip

In \S 5 we will see that the ideals of generic subspace
arrangements often fail to be pl-generated.
We study arrangements of points in ${\bf P}^2$ in \S \ref{sec:
points} and
   arrangements of lines in ${\bf P}^3$ in \S \ref{sec: lines}.
Propositions \ref{prop: points} and \ref{prop: lines} show that the ideals of
  generic arrangements of points in ${\bf P}^2$ (respectively lines in 
${\bf P}^3$) are
not  pl-generated when the number of subspaces is large.  However, in 
both cases
generic arrangements are scheme-theoretically cut
out by products of linear forms.  This is true for any arrangement of pairwise
disjoint subspaces; it is easy to see that the union of any
disjoint subschemes of projective space cut out by ideals $I_1,
\ldots, I_r$ is scheme-theoretically defined by $I_1 \cdots I_r.$  By contrast,
   Proposition \ref{prop: cones} shows that there exist arrangements of lines in
${\bf P}^3$ that are not scheme-theoretically cut out by any pl-generated
ideal.

%\vskip .5cm
\medskip

\noindent{\bf Acknowledgments.} We thank H. Derksen, D. Eisenbud, T. Geramita,
   D. Kozlov, R. Lazarsfeld and B. Sturmfels for helpful
conversations.  We also thank Adam Van Tuyl for pointing out a mistake in an earlier version.  

Anders Bj\"orner is supported by the G\"oran
Gustafsson Foundation for Research in
Natural Sciences and Medicine, and by MSRI, Berkeley. Irena Peeva is partially
supported by NSF and by MSRI, Berkeley.  Jessica Sidman is supported
  by an NSF postdoctoral fellowship and by MSRI, Berkeley.

%\vglue .8cm

\section{Notation and conventions}  

\noindent
We begin by briefly recalling a few basic definitions.
If $x, y$ are elements
of a geometric lattice $L$ then
$x \wedge y$ is their ${\it meet},$ or greatest lower bound in $L$, and $x \vee
y$ is their ${\it join},$ or least upper bound in $L.$
The least element of $L$ is denoted by $\hat{0}$,
and the greatest by $\hat{1}$.
A set $A \subseteq L$ is an ${\it antichain}$ if
$\hat{0}\notin A \neq \emptyset$ and the elements of $A$
are pairwise incomparable with respect to the partial ordering in $L.$

We say that a subspace arrangement $\AAA$ is {\em embedded}
in a hyperplane arrangement $\HHH$
if each $X\in\AAA$ is the
intersection of some of the hyperplanes from $\HHH$.
%cf. \cite[p. 328]{bjorner92}. 
Throughout the paper,
   ${\mathcal{A}}\subseteq {\mathcal{H}}$ is an embedding into a central
hyperplane arrangement $\mathcal{H}$
with intersection lattice $L_{\mathcal{H}}$. Denote by $\ell_1,\dots ,\ell_p$
the linear forms defining the
hyperplanes in $\mathcal{H}$. We think of $L_{\mathcal{H}}$ as the geometric
lattice with
atoms $\ell_1,\dots ,\ell_p$.
Denote by $V_1,  \dots , V_r$  the elements in
$L_{\mathcal{H}}$ that correspond to the subspaces in $\mathcal{A}$. The set
$\AAH
=\{ V_1,  \dots , V_r\}$
is an antichain in the lattice $L_{\mathcal{H}}$.

For a comprehensive
introduction to general notions related to hyperplane arrangements and subspace
arrangements see \cite{orlik-terao} and \cite{bjorner92}, respectively.  For
the matroid and geometric lattice
point of view see \cite{oxley}.

For simplicity, we assume that the field $k$ is infinite.  If $k$ is finite
then the ideal $I_{\mathcal A}$ should be defined as the intersection
of ideals $I_i$ generated by linear forms such that the Krull dimension of
$k[x_1,\ldots, x_n]/I_i$ is the dimension of the corresponding vector subspace
of $k^n.$  However, this may be strictly contained in the ideal of
all polynomials which vanish on the finitely many points of the subspaces.
%\vfil\break
\section{Blocker duality}

In this section we define the blocker dual of a subspace arrangement
embedded in a hyperplane arrangement and discuss properties of the associated 
blocker ideal. We begin and end the section by describing examples for which
the blocker ideal is the radical ideal of an arrangement.

\subsection{Motivating examples}\label{sec: motivation}

The main motivation for the construction which we will describe in \S
\ref{sec: blockers} comes from observing a beautiful duality between the
subspaces of certain arrangements embedded in the 
braid arrangement and polynomials which generate their
defining ideals. In order to discuss these examples we recall some
   basic facts about braid arrangements.  (See \cite{orlik-terao} and
\cite{bjorner92} for more details.)

The braid arrangement $\mathcal{H}_n$ is the arrangement of hyperplanes in
$k^n$ defined by the polynomial \[\prod_{1 \leq i < j\leq n} (x_i -x_j) .\] We
identify the intersection lattice of $\mathcal{H}_n$ with $\Pi_n,$
the lattice of
all partitions of $[n] = \{1, \ldots, n\}$, as follows.  Given a
partition $\pi$
   of $[n]$ into disjoint blocks, we define $i \equiv j$ if and only if $i$ and
$j$ are in the same block of $\pi$ and associate to $\pi$ the linear subspace
of $k^n$ defined by the ideal $(x_i -x_j \, |\, 1 \leq i < j \leq n,
\, i \equiv j).$

The symmetric group $S_n$ acts on the intersection lattice of $\mathcal{H}_n$
   by permuting the subscripts of the coordinates of $k^n.$  The orbits of this
action are indexed by the shapes of partitions of the set $[n].$  We
say that the
\emph{shape} of a partition $\pi$ is the list of its block sizes arranged
in non-increasing order.  E.g., if $\pi$ is the partition $\{
\{1,2,3\}, \{4\}\}$ of $[4],$
then the shape of $\pi$ is $(3,1).$

Let  $\mathcal{A}_{\lambda}$ be the arrangement consisting of all
subspaces corresponding to partitions
of shape $\lambda$.
The products \[f_{\pi} = \prod_{\stackrel{ i <j}{i \equiv j}}(x_i-x_j)\]
for $\pi\in\Pi_n$ play an
important role in what follows.

\begin{ex}\label{ex: li-li}
In \cite{li-li1} Li and Li found an explicit system of generators for 
the ideals
of orbit arrangements $\mathcal{A}_{\lambda}$ corresponding to ``hook'' shapes
   $\lambda = (m,1, \ldots, 1)$.
Namely, the vanishing ideal of the arrangement $\mathcal{A}_{\lambda}$ is
\[ (f_{\pi} \:| \:\pi\, \mathrm{has} \, m-1 \,\mathrm{blocks}).\]
\end{ex}

\begin{ex}\label{ex: lovasz}
A result of Kleitman and Lov\'asz in \cite{lovasz} describes a system of
generators of the ideal of certain arrangements consisting of unions of orbit
arrangements. Let \[\mathcal{A}_m = \bigcup_{\lambda \, \mathrm{has} \, m-1 \,
\mathrm{blocks}} \mathcal{A}_{\lambda}.\]  The defining ideal of
$\mathcal{A}_m$ is \[ \bigl(f_{\pi} \: | \:\pi \, \mathrm{has \, shape}\,
(m, 1, \ldots, 1) \bigr).\]
\end{ex}

Note that the partitions indexing the subspaces of Example \ref{ex:
li-li} index
the generators of the ideal of the arrangement in Example \ref{ex: lovasz}, and
vice versa.  In the next section
we define a combinatorial operation on antichains in a geometric
lattice which captures this duality.\medskip

\subsection{The blocker}\label{sec: blockers}
We now define the notion of \emph{blocker duality} motivated by
Examples \ref{ex: li-li} and \ref{ex: lovasz} and demonstrate its basic
properties. This concept is purely combinatorial and 
for our purposes best discussed
in the setting of geometric lattices. For more about the
combinatorial properties of the blocker construction, see
\cite{matveev} and \cite{bjorner-hultman}.

\begin{defin}
Let $A$ be an antichain in a geometric lattice $L$. The \emph{
blocker} of $A$ is the antichain
$$
A^*= \mathrm{\, min \,}  \, \{\, x\in L
\, | \,
a\wedge x \neq \hat{0} \,
\mathrm{ for \, every }\,  a\in A\,\}\, ,
$$
where $\mathrm{\, min \,} E$ denotes the set of minimal
elements of a subset $E\subseteq L$.
\end{defin}

Note that $A^*\neq\emptyset$, since $a\wedge \hat{1}=a\neq \hat{0}$
for all $a\in A$.
As an example, let $A=\{\hat{1}\}$ and $B=\{\mathrm{atoms}\}$.
Then $A^*=B$ and $B^*=A$.

A partial order
on the antichains in a geometric lattice $L$ is
defined as follows:
we say that $A\le B$ for two antichains if for each $b\in B$ there
exists an $a\in A$ such that
$a\le b$.
The proof of the following lemma is straightforward.

\begin{lem}\label{lem: duality}
\begin{enumerate}
\item[]
\item[(1)] If $A\le B$, then $B^*\le A^*$.
\item[(2)] $A^{**}\le A$.
\end{enumerate}
\end{lem}

The following proposition describes the sense in which the $*$ operation
on antichains is a reflexive duality operation.
Note that the notion of reflexivity given by $*$ is somewhat weak: $A = A^{**}$
   does not hold in general.

\begin{prop}\label{prop: duality}
Let $A$ be an antichain in a
   geometric
lattice $L$. Then $A^*=A^{***}$.
\end{prop}

\begin{proof}
By Lemma \ref{lem: duality} (2)
we get that
$A^{***}\le A^*$. On the other hand, Lemma \ref{lem: duality} (1) applied to
$A^{**}\le A$ yields
$A^{***}\ge A^*$.
\end{proof}

The definition of the blocker is designed to generalize the duality between
Examples \ref{ex: li-li} and \ref{ex: lovasz}. Indeed, we have:

\begin{ex}\label{ex: braid duals}
Let $\lambda = (m, 1, \ldots, 1)$, and let
   $A_{\lambda}$ be the antichain in $\Pi_n$
of partitions of the set $[n]$ of shape $\lambda$.
Then $$A^*_{\lambda}=\{\pi\in\Pi_n\,\mid\, \pi \;\mathrm{has} \; m-1 \;
\mathrm{blocks}\}
$$
and ${A^{**}_{\lambda}} = A_{\lambda}.$
The antichain $A_{\lambda}$
corresponds to
the orbit arrangement $\mathcal{A}_{\lambda}$ embedded in
   the braid arrangement  $\HHH_n$, and $A^*_{\lambda}$ corresponds to
   the arrangement \[\mathcal{A}_m = \bigcup_{\lambda \, \mathrm{has} \, m-1 \,
\mathrm{blocks}} \mathcal{A}_{\lambda}\, .\]
\end{ex}

\begin{rmk}\label{ex: bool duals}
The blocker construction was originally introduced (in 
\cite{edmonds-fulkerson} and other places) for the special case when
$L$ is the Boolean lattice of all subsets of a finite set $V$.
%The blocker $B$ of an antichain $A$ in $2^V$ is characterized by the property:
%{\em For all $U\subseteq V$, $U$ contains a member of $A$
%if and only if $V\setminus U$ does not contain a member of $B$.}
In this case it is known that $A^{**}=A$ for all antichains $A$.
See Example \ref{ex: coord blocker} for more about this.

The generalization of blockers to posets has also independently
been considered by Matveev \cite{matveev}.
\end{rmk}
\vglue .4cm

\subsection{The blocker ideal}\label{sec: blocker ideal}

We now define the blocker ideal $\BAxH=\BAH$ of an arrangement
$\AAA$, with respect to an embedding
   $\mathcal{A} \subseteq \mathcal{H}$, and show
some of its most basic connections with
the vanishing ideal $\IA$.

Suppose that a subspace arrangement $\AAA$ is embedded
in a hyperplane arrangement $\HHH$ , i.e. $X\in\LH$ for all $X\in\AAA$.
Let $\AAH$ denote $\AAA$ viewed as an antichain in $\LH$
and $\AaAH$ denote its blocker dual. When confusion cannot arise we simplify
notation by suppressing the reference to $\HHH$ and identifying
the subspace arrangements embedded in $\HHH$ with the antichains
contained in $\LH$. Thus, we may speak directly of the blocker dual
$\AAA^*$ of a subspace arrangement $\AAA$.
Note that the operation $\AAA\rightarrow\AAA^{**}$
defines a closure operation on subspace arrangements
(with respect to $\HHH$), namely,
by Lemma~\ref{lem: duality}: $V_{\AAA}\subseteq V_{\AAA^{**}}$.

Given an arrangement of hyperplanes $\mathcal{H}$ in which the hyperplanes are
defined by linear forms $\ell_1 , \ldots, \ell_p,$ we associate a
product of linear
forms to each element of of $L_{\mathcal{H}}$ as follows:

\begin{defin}\label{defin: product}
For $X \in L_{\mathcal{H}},$ define \[Q_{X} = \prod_{\ell_i(X) \equiv 0}
   \ell_i.\]
\end{defin}

Using this definition
we  define the blocker ideal of
$\mathcal{A} \subseteq \mathcal{H}:$

\begin{defin}\label{defin: blocker ideal}
The \emph{blocker ideal} $\BAH$
is   \[\BAH = ( Q_X \, | \, X \in \AAA^*).\]
\end{defin}
\noindent
The first part of the following proposition
shows that  that
$Q_{X} \in I_{\mathcal{A}}$
for all $X \in \AAA^*,$
hence
$$\BAH \subseteq I_{\mathcal{A}}.$$
The third part shows that
the blocker ideal cuts out $V_{ \AAA^{**}  }$
set-theoretically.

Let $\AAA$ be a subspace arrangement embedded in the hyperplane arrangement
$\HHH$, and let $\AAA^*$ be the blocker dual arrangement. For $X\in \LH$ let
$\HHH /X = \{H\in \HHH \mid H\supseteq X\}$.

\begin{prop}\label{settheor}
\begin{enumerate}
\item[]
\item[(1)]  $V_{\AAA} \,\,\,\,\,\subseteq \,\bigcap_{X\in \AAA^*} V_{\HHH /X}$
\item[(2)]  $V_{\AAA^*} \,\,\,\supseteq \,\bigcap_{X\in \AAA} V_{\HHH /X}$
\item[(3)]   $V_{\AAA^{**}} \, = \,\bigcap_{X\in \AAA^*} V_{\HHH /X}$
\end{enumerate}
\end{prop}

\begin{proof}
(1) Suppose that $z\in C\in \AAA$. For each $X\in \AAA^*$
there exists (by definition of the blocker) a hyperplane $H_X\in \HHH$
such that $z\in C\cup X\subseteq H_X$. Then,
   $z\in \cap_{X\in \AAA^*} H_X\,\,\subseteq \,\cap_{X\in \AAA^*} V_{\HHH /X}.$

(2) Suppose that $y\in V_{\HHH /X}$ for all $X\in \AAA$.
%$and $x\neq 0$.
So, for each $X\in \AAA$ there is a hyperplane $H_X \supseteq X$ such that
$y\in H_X$. Let $C=\cap_{X\in \AAA} H_X$. Then $C\wedge X \neq \hat{0}$
for all $X\in \AAA$, and hence there exists some
$\widetilde{C}\in\AAA^*$ such that
$\widetilde{C} \le C$. We have that $y\in C\subseteq
\tilde{C}\subseteq V_{\AAA^*}$.

(3) Using the preceding parts we have that
\begin{eqnarray*}
V_{\AAA^{**}} \, \subseteq\,  \bigcap_{X\in \AAA^{***}} V_{\HHH /X}
   \,=\,\bigcap_{X\in \AAA^*} V_{\HHH /X}\, \subseteq\, V_{\AAA^{**}}.
\end{eqnarray*}
\end{proof}

\begin{thm}\label{thm: blocker ideal 1} 
Over an algebraically closed field, the following properties hold:
\begin{enumerate}
\item[(1)]  $\rad (\BAH) = I_{\mathcal{A^{**}}},$
\item[(2)]  $\rad (\BAH) = I_{\mathcal{A}} \mbox{ \;  if and only if \;      }
\AAA^{**} = \AAA.$
\end{enumerate}
\label{thm: A starstar}
\end{thm}

\begin{proof}
The first statement follows directly from Proposition
\ref{settheor}(3) via the Hilbert Nullstellensatz.
Thus,   $\rad (\BAH)$ defines $V_{\AAA^{**}}$ and so  the second part
follows.
   \end{proof}

As we will see with Example \ref{ex: kozlov}, the property
$\AAA =\AAA^{**}$ does not guarantee that
$\BAH$ is a radical ideal.
\medskip

\subsection{More examples}\label{sec: examples}
The notion of blocker duality  behaves well for several
interesting subspace arrangements.
Here we give examples having the property
that $\BAH = I_{\mathcal{A}}$, assuming only that the field $k$ is infinite.

\begin{ex}\label{ex: hyperplane blocker}
Suppose that $\mathcal{H} = \{H_1, \ldots, H_p\}$ is a central hyperplane
arrangement with defining equation $\ell_1 \cdots \ell_p.$  Then $X:=
\cap_{i = 1}^{p} H_i$ is an element of $L_{\mathcal{H}}$ and is the only
element of $\HHH^*.$ Furthermore,
$Q_X = \ell_1 \cdots \ell_p$. Hence,
$B_{\mathcal{H}, \mathcal{H}}= I_{\mathcal{H}}$.
\end{ex}

\begin{ex}\label{ex: coord blocker}
We say that $\mathcal{A}$ is a {\em coordinate subspace arrangement}, or a
\emph{Boolean} arrangement (see \cite[\S 3.2]{bjorner92}), in $k^n$ if each
subspace in $\mathcal{A}$ is an intersection of coordinate hyperplanes.
Such an $\AAA$ has a natural embedding into the coordinate hyperplane
\newcommand{\CCC}{\mathcal{C}_n}
arrangement $\CCC$ defined by the ideal $(x_1 \cdots x_n)$, whose
intersection lattice
is isomorphic to the Boolean lattice $B_n$ of all subsets of $[n]$.
The blocker duals of coordinate subspace arrangements
have close connections with the Stanley-Reisner rings of simplicial complexes
and a nice interpretation in terms of Alexander duality, as
we now show.

Let us begin set-theoretically. An antichain $A$ in $B_n$ generates an abstract
simplicial complex $\Delta_A =\{X\subseteq [n] \mid X\subseteq F
\mbox{ for some }
F\in A\}$. Conversely,
$\mathrm{\, max \,} (\Delta)$ is an antichain for every simplicial complex
$\Delta$. Clearly, $$\Delta_{\mathrm{\, max \,} (\Delta)}=\Delta
\mbox{\quad and\quad}
\mathrm{\, max \,} (\Delta_A)=A,$$
so antichains and simplicial complexes are interchangeable concepts here.

Let $X^c =[n]\setminus X$ for subsets $X\subseteq [n]$,
and $A^c = \{X^c \mid X\in A\}$ for antichains $A$.
The simplicial complex
$\Delta^{\mathrm{dual}}=\{X^c \mid X\notin \Delta\}$ is known as the
{\em (combinatorial) Alexander dual} of $\Delta$. By the previous
comments we may instead speak of the Alexander dual $A^{\mathrm{dual}}$
of an antichain $A$ in $B_n$.

We know from Example \ref{ex: bool duals}
that $A^{**}=A$ for all antichains $A$ in $B_n$.
We also know that $(A^{\mathrm{dual}})^{\mathrm{dual}}=A$
and $A^{c\, c}=A$ for all $A$. These duality
operations are related as follows
\begin{equation}
A^{\mathrm{dual}}=A^{c\, *\, c}  \mbox{\quad and\quad}
A^* = A^{c\, (\mathrm{dual})\, c},
\end{equation}
since the definitions show that
\begin{equation}\label{duals}
X\in A^{ (\mathrm{dual})\, c} \;
\Leftrightarrow \; X\in \mathrm{\, min \,} (B_n\setminus \Delta_A)  \;
\Leftrightarrow  \;   X\in A^{c\, *}.
\end{equation}

Let
$$G\subseteq [n] \quad \leftrightarrow \quad
S_G =\{(x_1, \dots , x_n)\in k^n\mid x_i =0 \mbox{  for all  } i\notin G\}$$
be the chosen correspondence between subsets of $[n]$ and coordinate
subspaces.  If $\mathcal{A}$ is a coordinate subspace arrangement 
corresponding to an
antichain $A$ in $B_n$ with
simplicial complex $\Delta_A ,$ then $I_{\mathcal{A}}$ is a monomial ideal.
Namely, the ideal $I_{\mathcal{A}}$ is generated by the square-free products of
variables $\prod_{i \in G} x_i,$ for all minimal non-faces $G\notin\Delta_A$
(see \cite[\S 11.1]{bjorner92}).
This is known as the {\em Stanley-Reisner ideal} of  $\Delta_A$.

Under the order-reversing isomorphism $B_n \leftrightarrow L_{\CCC}$
given by $G\leftrightarrow S_G$ we have that $A\leftrightarrow \AAA$ implies
that $A^{c\, *\, c}\leftrightarrow \AAA^*$. Furthermore,
$$\prod_{i \in G} x_i =Q_{S_{G^c}}$$
for all $G\subseteq [n]$.
Equation (\ref{duals}) shows that
$A^{c\, *} =\mathrm{\, min \,} (B_n\setminus \Delta_{A}),$ from
which follows that $\IA$ is generated by all $Q_{S_G}$
such that $G\in A^{c\, *\, c}$. That is, $\IA$ is in fact the blocker ideal.

Hence, for all coordinate subspace arrangements:
$B_{\AAA, \CCC} = I_{\mathcal{A}}$.
\end{ex}

\begin{ex}The results of \cite{li-li1} and \cite{lovasz} show that the blocker
ideals of $\mathcal{A}_{\lambda}$,
for $\lambda$ of hook shape, and of its *-dual in the braid arrangement, are
the respective radical ideals, cf.
Examples \ref{ex: li-li}, \ref{ex: lovasz} and \ref{ex: braid duals}. Orbit
arrangements $\mathcal{A}_{\lambda}$ are themselves in general not blockers with
respect to the braid arrangement. 
A procedure for computing their
blocker duals, and hence their blocker ideals,
is given in \cite{bjorner-hultman}. 
We do not know of any description
of their vanishing ideals for general non-hook shapes

The following table gives the blocker duals, and double duals,
with respect to the braid arrangement, for all
$\mathcal{A}_{\lambda}$ indexed by partitions $\lambda$ of $n=6$
that are not of hook shape.
Here $\ll m\gg$ denotes the union of all orbit arrangements for partitions with
$m$ blocks, as in Example \ref{ex: lovasz}.
\medskip

\begin{center}
$\begin{array}{c|c|c}
\mathcal{A}_{\lambda} & \mathcal{A}_{\lambda}^* & \mathcal{A}_{\lambda}^{**}\\
\hline
(4,2) & (2,2,2) \cup (3,1,1,1) & (4,2) \cup (5,1)\\
(3,3) & (3,1,1,1) & \ll 2\gg\\
(3,2,1) & (3,3) \cup (4,1,1) & (3,2,1) \cup (4,1,1)\\
(2,2,2) & (4,1,1) & \ll 3\gg\\
(2,2,1,1) & (5,1) & \ll 4\gg
\end{array}$
\end{center}
\end{ex}

Note that the arrangements
$\AAA_{(2,2,2)} \cup \AAA_{(3,1,1,1)}$ and $\AAA_{(4,2)} \cup \AAA_{(5,1)}$
are blocker dual to each other, as are the arrangements
$\AAA_{(3,3)} \cup \AAA_{(4,1,1)}$ and $\AAA_{(3,2,1)} \cup \AAA_{(4,1,1)}.$
Using MACAULAY 2 \cite{grayson-stillman} we computed the ideals of
these four arrangements
and compared them to the respective blocker ideals. Working over the field
$\mathbb{Q}$ we found that the blocker ideal in each case
equals the vanishing ideal.

\begin{ex}
Let $E$ be a $d$-dimensional vector space over the field $k$. Given
positive integers
$m$ and $n$ and a function $f: [m] \rightarrow [n]$ let
$$W_f = \{ ( x_1 ,\dots , x_n, x_{f(1)}, \dots , x_{f(m)})\mid x_i \in E\}.
$$
This is an $nd$-dimensional linear subspace of $E^{n+m}$. Letting $f$
range over all such functions, define the {\em polygraph arrangement}
\newcommand{\ZZZ}{\mathcal{Z}}
$$\ZZZ_E (n,m) = \{W_f \mid f: [m]\rightarrow[n]\}.
$$
Such arrangements were introduced by
M. Haiman in \cite{haiman}, and for $E=\mathbb{C}^2$
they play a crucial role in his proof of the $n!$ conjecture. They were further
investigated from a combinatorial point of view in \cite{hultman}.

Now let $d=1$, and consider the vanishing ideal $I_{\ZZZ_k (n,m)}$
in the polynomial ring $k[x_1,\dots,x_n,a_1, \dots, a_m]$. Haiman 
\cite[p. 966]{haiman}
shows that the ideal $I_{\ZZZ_k (n,m)}$ is generated by 
$$q_i=\prod_{j\in [n]} (x_j -a_i), \quad i\in [m].
$$
This implies that $I_{\ZZZ_k (n,m)}$ is the blocker ideal of $\ZZZ_k
(n,m)$ with respect to its
embedding into the ``bipartite braid arrangement''
$\HHH (n,m)= \{x_j -a_i\mid j\in [n], i\in [m]\}$, as we now show.

Let $P$ and $Q$ be disjoint sets of cardinalities $|P|=m$ and $|Q|=n$, and let
$\Pi_{P\cup Q}$ denote the lattice of all partitions of the set $P\cup Q$.
\newcommand{\PPP}{\Pi^{\circ}_{P,Q}}
Let $\PPP$ denote the lattice that is join-generated within $\Pi_{P\cup Q}$
by all rank one partitions (atoms) whose only non-singleton block
is of type $\{p,q\}$ with $p\in P$ and $q\in Q$.
So, $\pi\in\PPP$ if and only if every block of  $\pi$ either is a
singleton or else
intersects both $P$ and $Q$.

The isomorphism of the intersection lattice of the braid arrangement
$\HHH_{n+m}$
with $\Pi_{P\cup Q}$ (see \S 3.1) clearly restricts to an isomorphism
$L_{\HHH (n,m)} \cong \PPP$. Hence, we can compute the blocker dual
of $\ZZZ_k (n,m)$ within $\PPP$.

The antichain in $\PPP$ that corresponds to $\ZZZ_k (n,m)$ under the stated
isomorphism is the antichain $A_{n,m}$ of all partitions in $\PPP$ for which
each block contains {\em exactly one} element from $Q$ (note that
consequently  each $\pi\in A_{n,m}$
contains exactly $n$ blocks). Via combinatorial reasoning one sees that
\begin{eqnarray*}
A^*_{n,m} \!\!\! &=& \!\!\! \{\pi\in\PPP\mid \pi \mbox{ has a unique
non-singleton block }
Q\cup\{p_i\},\\
& & \qquad\qquad\qquad  \mbox{ for some } i\in [m]\}\\
A^{**}_{n,m} \!\!\! &=&\!\!\! A_{n,m}.
\end{eqnarray*}
Hence, the generator set of the blocker ideal $\BAH$, for
$\AAA=\ZZZ_k (n,m)$ embedded in  $\HHH=\HHH (n,m)$, is precisely
the set of polynomials $\{q_i\}_{i\in [m]}$ defined above.
In other words, $\BAH=\IA$ for polygraph arrangements
in case $d=1$.

The situation becomes more complicated if $d>1$. The combinatorics stays
the same, but the algebra gets more involved. Haiman \cite[\S 4.6,
eq. (96)]{haiman}
gives a set of generators for the special case of $n=d=2$,
but states \cite[p. 967]{haiman} that  ``at present, we do not have a good
conjecture as to a set of generators for the full ideal [for $d=2$]
in general''.
\end{ex}

\vglue .7cm

\section{Products of linear forms inside $I_{\AAA}$}
In this section we show that the ideal generated by all products
of linear forms that vanish on an arrangement $\AAA$ can be constructed
by a simple algorithmic procedure.  Given an embedding $\AAA \subset \HHH,$ 
we construct the ideal $F_{\AAA, \HHH}$ generated by all products of linear forms 
defining elements of $\HHH.$   We then show how to generate an embedding 
of $\AAA$ into a hyperplane arrangement so that the ideal $F_{\AAA,\HHH}$ is 
as large as possible.

\subsection{The  $\HHH$-product ideal}
Example \ref{ex: blocker failure} shows that the blocker ideal
$B_{\AAA,\HHH}$ may fail to be equal to $I_{\AAA},$ as the blocker
construction may not
detect all products of linear forms in $I_{\AAA}.$  We introduce the {\it
$\mathcal{H}$-product ideal} $F_{\mathcal{A},\mathcal{H}}$,
which corrects for this failure.

\begin{defin}\label{def: full product}
The {\it  $\mathcal{H}$-product ideal} is the ideal
\begin{eqnarray*}
F_{\mathcal{A},\mathcal{H}}=
\biggl({\ell}_1\cdots {\ell}_q\,  & \mid & {\ell}_j \in L_{\HHH}, \;
\hbox{and for all $1\le i\le r$} \\ & &
\hbox{there exists $1\leq j \leq q$
such that ${\ell}_{j}\in I_i$}
\ \,\biggr)\, .
\end{eqnarray*}
\end{defin}

Clearly, $B_{\mathcal{A},\mathcal{H}}\subseteq
F_{\mathcal{A},\mathcal{H}}\subseteq \IA,$ and the first two
ideals can be  computed combinatorially given $L_{\mathcal{H}}$
and the antichain $\AAH$.
We will  show that a strict inclusion $B_{\mathcal{A},\mathcal{H}}\subset
F_{\mathcal{A},\mathcal{H}}$ is possible.
That strict inclusion
$F_{\mathcal{A},\mathcal{H}}\subset \IA$
is possible can be seen from Example \ref{ex: kozlov}.

\begin{ex} \label{ex: blocker failure}
Let
$\mathcal{H}=\HHH_3$ be the braid arrangement with hyperplanes
$\{x_1=x_2\},\ \{x_1=x_3\},\ \{x_2=x_3\}$.
Consider the subspace arrangement $\mathcal{A}$ with subspaces
$\{x_1=x_2\},\ \{x_1=x_3\}$. Clearly,
${\mathcal{A}}\subseteq{\mathcal{H}}$.
   Both ideals
$B_{\mathcal{A},\mathcal{H}}$ and $F_{\mathcal{A},\mathcal{H}}$ are principal,
but they are different, and furthermore their radicals are different as well.
We have that
$$
B_{\mathcal{A},\mathcal{H}}=\biggl((x_1-x_2)(x_1-x_3)(x_2-x_3)\biggr)
\subset
\biggl((x_1-x_2)(x_1-x_3)\biggr)=F_{\mathcal{A},\mathcal{H}}= I_{\AAA} .$$
\end{ex}

The ideal $F_{\mathcal{A},\mathcal{H}}$
can be characterized algebraically as follows.

\begin{prop} \label{prop: full product ideal}
   The ideal  generated by all products of
linear forms in
$L_{\mathcal{H}}$ that are contained in $I_{\mathcal{A}}$ is equal to
$F_{\mathcal{A},\mathcal{H}}.$
\end{prop}

\begin{proof}
   Consider
a product of linear forms ${\ell}_1\cdots {\ell}_q$ such that each ${\ell}_i\in
L_{\mathcal{H}}$. We have that
${\ell}_1\cdots {\ell}_q\in I_{\mathcal{A}}$, if and only if for each
$1\le i\le r$
there exists a linear form
${\ell}_{j}\in I_i.$
\end{proof}

\vglue .2cm

\subsection{The embedding}
The next result shows that the vanishing ideal of every arrangement of
   two subspaces is pl-generated. This is not true for  three
subspaces, as shown by an example due to Li and Li \cite{li-li2}; see
Proposition \ref{prop: lines}  and Remark \ref{rmk: lines} for
comments and generalizations.

\begin{prop}\label{prop: 2 subspaces} 
Let $J$ and $J'$ be two linear ideals. The
ideal  $J\cap J'$ is generated by products of linear forms.
\end{prop}

\begin{proof}
%A linear change of coordinates on the ambient space takes $J$ and $J'$ to 
%monomial ideals and the result follows.
Let $V=k^n$ denote the ambient vector space, and let $V_1$ and $V_2$ be the
vector subspaces of $V$ defined by $J$ and $J',$
respectively.  Write $V_1 = W_1 \oplus (V_1 \cap V_2),$
$V_2 = W_2 \oplus (V_1 \cap V_2),$
and $V = W \oplus W_1 \oplus W_2 \oplus (V_1 \cap V_2),$
for vector spaces $W_1 \subseteq V_1,$ $W_2 \subseteq V_2,$ and $W 
\subseteq V.$
Let $\{s_1,\dots, s_q\}$ be a basis for $W,$ $\{e_1,\dots ,e_b\}$ be a basis
for $W_1,$ $\{h_1,\dots,h_a\}$ be a basis for $W_2,$ and $\{ t_1, \dots, t_c\}$
be a basis for $V_1 \cap V_2.$

Clearly, $J = (s_1,\dots, s_q,  h_1,\dots,h_a)$ and
  $J' = (s_1,\dots, s_q, e_1,\dots ,e_b)$.  Since $J$ and $J'$ are monomial
ideals, their intersection is generated by $s_1, \dots, s_q$
  and all elements of the form $e_ih_j$.
\end{proof}

However, even for an arrangement of two subspaces, one can choose a
poor embedding into a
hyperplane arrangement from which one cannot readily
detect if $I_{\mathcal{A}}$ is generated by products of
linear forms:

\begin{ex}[D. Kozlov]\label{ex: kozlov} 
Suppose that the characteristic  of $k$ is not equal to 2.
Consider the subspace arrangement ${\mathcal{A}}$ that consists of  the two
subspaces
$\{ x_1-x_2=0,\  x_1+x_2=0\}$ and $ \{ x_1-x_3=0,\  x_1+x_3=0\}$.
Take the hyperplane arrangement
$\mathcal{H}$ consisting of the hyperplanes $\{ x_1-x_2=0\} , \ \{ x_1+x_2=0\}
,\
    \{ x_1-x_3=0\} $ and $\{ x_1+x_3=0\}$. Then
$ A_{\mathcal{A},\mathcal{H}}=A_{\mathcal{A},\mathcal{H}}^{**}$. On
the other hand,
\begin{eqnarray*}
   B_{\mathcal{A},\mathcal{H}}=F_{\mathcal{A},\mathcal{H}}
&=\Bigl((x_1-x_2)( x_1-x_3),\  (x_1-x_2)( x_1+x_3), \\
            &\ \ \ \ \ (x_1+x_2)(x_1-x_3),\  (x_1+x_2)(x_1+x_3) \, \Bigr)\\
   & =\Bigl( x_1^2,x_1x_2,x_1x_3,x_2x_3 \, \Bigr)
\end{eqnarray*}
is clearly not a reduced ideal, so it is not equal to $I_{\mathcal{A}}$.

By contrast, let us consider the embedding of $\mathcal{A}$ into the
coordinate hyperplane arrangement
$\mathcal{C}_3$ consisting of the hyperplanes $\{ x_1=0\},\ \{ x_2=0\},\
\{ x_3=0\}$.
Denote by $J$ and $J'$ the defining ideals of the two subspaces in
$\mathcal{A}$.
In this case, $\{ x_1\} $ is a basis of the $k$-space $J_1\cap J_1'$.
Furthermore, $\{ x_1,x_2\}$ is a basis of the $k$-space $J_1$,
    and $\{ x_1,x_3\} $ is a basis of the $k$-space $J_1'$. By
Proposition \ref{prop: 2 subspaces},
we conclude that
$$I_{\mathcal{A}}=(\, x_1,\ x_2x_3)=B_{\mathcal{A},\mathcal{C}_3}\,
.$$
This can also be seen to follow from Example \ref{ex: coord blocker}.
\end{ex}

The problem with the first embedding in Example \ref{ex: kozlov} is that the
ideals $B_{\mathcal{A},\mathcal{H}}$  and $F_{\mathcal{A},\mathcal{H}}$ are
strictly smaller than the ideal  generated by all products of linear
forms in $I_{\mathcal{A}}$. We will show that this problem can be
avoided if we
take an embedding into a larger
hyperplane arrangement ${\tilde{\mathcal{H}}}\supseteq {\mathcal{H}}$.

\begin{con}\label{con: embed}
We start with the embedding $\mathcal{A}\subseteq \mathcal{H}$. If necessary,
to the atoms
$\ell_1,\dots ,\ell_p$ of $L_{\mathcal{H}}$ we add finitely many new atoms
to obtain a new larger
hyperplane arrangement $\tilde{\mathcal{H}}$ such that for any
choice of $1\le i_1 <\dots < i_t\le r$ there exists
a subset of atoms in $L_{\tilde{\mathcal{H}}}$ that forms a basis for  the
$k$-space
consisting of the linear forms in $I_{i_1}\cap\dots \cap I_{i_t}$.
\end{con}

The procedure for enlarging $\HHH$ to $\tilde\HHH$ is clearly finite,
since it amounts to adding a finite number of atoms (linear forms)
in each of at most $2^r$ steps. The key observation is that the ideal
$F_{\AAA,\tilde\HHH}$ is the largest possible, and is independent 
of the choice of $\HHH$ and $\tilde\HHH$, as we now show.

\begin{thm}\label{prop: full ideal}
The ideal  generated by all products of linear forms in $I_{\mathcal{A}}$ is
equal to the  $\tilde{\mathcal{H}}$-product ideal
$F_{\AAA,\tilde\HHH}$.
\end{thm}

\begin{proof}
A product $f=f_1\cdots f_q$ of linear forms is in $I_{\mathcal{A}}$ if and
only if
for each $0\le i\le r$ there exists a linear form
    $f_{j}\in I_i$. Furthermore,
if $f_q\in I_{i_1}\cap \dots \cap I_{i_t}$, then $f_q$ is a linear
combination of the basis elements of the $k$-space consisting of all
linear forms in $I_{i_1}\cap \dots \cap I_{i_t}$. Hence,
the ideal  generated by all products of linear forms
in $I_{\mathcal{A}}$ is
$F_{\AAA,\tilde\HHH}$.
\end{proof}

B. Sturmfels asked if one can check whether $I_{\mathcal{A}}$ is generated by
products of linear forms algorithmically.  We obtain such an algorithm as an
immediate corollary of Theorem \ref{prop: full ideal}. It can be
implemented using the computer algebra system
MACAULAY~2 \cite{grayson-stillman}. Note that our algorithm avoids
    computing  radicals, which are very  difficult to compute.

\begin{alg}\label{alg}
A subspace arrangement $\mathcal{A}$ is
given.
\item{(1)} Compute  $I_{\mathcal{A}}$ as the intersection of the
linear defining ideals \\
\indent \; of the subspaces in $\mathcal{A}$.
\item{(2)} Choose
an embedding into a hyperplane arrangement $\mathcal{H}$.
\item{(3)} Construct $\tilde{\mathcal{H}}$.
\item{(4)} Construct the  $\tilde {\mathcal{H}}$-product ideal
$F_{\AAA,\tilde\HHH}$.
\item{(5)} Check if the ideals $I_{\mathcal{A}}$ and
$F_{\AAA,\tilde\HHH}$  have the same Hilbert function. \\
\indent \; If {\em YES}: $I_{\mathcal{A}}$ is generated by products of linear forms. \\
\indent \;  If {\em NO}: \:It is not.
\end{alg}
\bigskip

\section{Arrangements in  ${\bf P}^2$ and ${\bf P}^3$}

The results in this section show that in general, an arrangement of points in ${\bf P}^2$
or lines in ${\bf P}^3$  will not have a pl-generated ideal.  

\subsection{Points in ${\bf P}^2$}\label{sec: points}

Let $S = k[x,y,z].$   It is relatively easy to see by direct computation that if 
${\mathcal{A}}$ is any
set of $r$ points in  ${\bf P}^2$ with $r \leq 4$ then $I_{\mathcal{A}}$ is
pl-generated.  The possible configurations can be organized according to
the maximum number of collinear points.  We leave the computation aside.

What happens when $r > 4?$  Recall that a set $\AAA$ of $r$ points in ${\bf P}^2$ is 
\emph{linearly general}
if no three are collinear  and that a set of $r$ points is $\emph{generic}$
if $\dim_k (S/I_{\AAA})_t = \min \{r, {t+2 \choose 2} \}.$  Five points in  ${\bf P}^2$ in linearly
general position lie on a unique irreducible conic, so their ideal cannot be
pl-generated. However, 6 generic points in linearly general position do not lie on 
a conic, and
we show in Proposition \ref{prop: 6points} that the ideal of such an 
arrangement of points is pl-generated. 
For $r >6,$ Proposition
\ref{prop: points} shows  that  the ideal of $r$ 
linearly general points in ${\bf P}^2$ is not pl-generated.  

It would be interesting to find a 
characterization of all
sets of points in
${\bf P}^2$ whose ideals are pl-generated.  In Proposition \ref{prop: 
pointsonconic}
we give an example of a constraint that one can impose on the geometry
of the points that forces their ideal to be pl-generated.
\smallskip

We begin by recalling some information about the ideals of points in ${\bf P}^2.$
A good reference for these results, which we will cite without proof, 
is Chapter 3 of \cite{eisenbud}.

The ideal of a finite set of points $\AAA$ in
${\bf P}^2$ has a very beautiful description via the Hilbert-Burch Theorem.
 Let $S(-d)$ denote the polynomial ring $S$ with
degrees shifted so that it is generated in degree $d,$ i.e., the 
degree $m$ piece is
  $S(-d)_m=S_{m-d}$.
The Hilbert-Burch Theorem says that
the ideal $I_{\AAA}$ is minimally generated by a nonzerodivisor $\alpha$ times
the  maximal minors of a $(t+1) \times t$ matrix $M$
that can be viewed as a map in the following short exact sequence:
\[
\xymatrix{ 0 \ar[r] & \overset{t}{\underset{i = 1}\bigoplus} S(-b_i) \ar[r]^M &
  \overset{t+1}{\underset{i = 1}\bigoplus} S(-a_i) \ar[r] &
I_{\AAA} \ar[r]& 0,}\]  where the $a_i$ are the degrees of the 
elements in a minimal
system of  generators of $I_{\AAA}$ and the $b_i$ are the degrees of 
the elements in a
minimal system of generators of the syzygies on the generators of $I_{\AAA}.$

All of the numerical information associated to $I_{\AAA}$ is encoded in the
degrees of the entries along the two main diagonals of $M.$  Let $e_i$
denote the degree of the $(i,i)$ entry of $M$ and let $f_j$ denote
the degree of the $(j,j+1)$ entry of $M.$  The following
theorem collects some of the relationships between the numbers we have
defined (see Proposition 3.8 in \cite{eisenbud}, for a proof of (1), (2), and (3)). 

\begin{thm}\label{thm: invariants}
Assume that $a_1 \ge a_2 \ge \cdots \ge a_{t+1}$ and $b_1 \ge b_2 \ge \cdots
\ge b_t.$   The following properties hold:
\begin{enumerate}
\item[(1)]  $e_i, f_i \ge 1$
\item[(2)] $a_i = \sum_{j<i} e_j + \sum_{j \ge i} f_j$
\item[(3)]  $b_i = a_i +e_i$
\item[(4)] $\deg \AAA = \sum_{i \leq j} e_if_j$
  (Ciliberto-Geramita-Orecchia
\cite{ciliberto-geramita-orecchia}).
\end{enumerate}
\end{thm}

We will  need the following corollary of the Hilbert-Burch Theorem.

\begin{cor}[Burch]\label{cor: numgens}
If a finite set of points in ${\bf P}^2$ lies on a curve of degree $d$ then
the ideal of the points can be generated by $d+1$ elements.
\end{cor}

 Programs in MACAULAY 2 \cite{grayson-stillman}, one of which
  was written by D.
Eisenbud, suggested that the ideal of six randomly chosen
points in ${\bf P}^2$ is pl-generated, motivating the following theorem.

\begin{prop}\label{prop: 6points} If $\AAA$ is a set of 6
generic and
linearly general points in ${\bf P}^2,$ then $I_{\AAA}$ is pl-generated.
\end{prop}

\begin{proof}
First we will show that $I_{\AAA}$ must be generated by 4
linearly independent cubics.  Then we will construct 4 degree 3
products of linear forms that vanish on $\AAA$ and are linearly independent.

Since six generic points impose six independent conditions on cubics, and
the space of cubics in three variables has dimension 10, we see that there are
precisely 4 linearly independent cubics in $I_{\AAA}.$  If the points are
chosen generically, they will not all lie on a line or 
a conic.  Thus, there are no elements in $I_{\AAA}$ of degree $\leq 2.$ By Corollary \ref{cor: numgens}, $I_{\AAA}$ requires at most 4 generators.   
Therefore, we see that the 4 cubics in $I_{\AAA}$ generate the ideal.

We construct 4 degree 3 forms vanishing on $\AAA.$  Label the points $p_1, \ldots, p_6$ and let
$L_{i,j}$ denote the line joining $p_i$ to $p_j.$  Since the points are
linearly general, the set of all $L_{i,j}$ with $i <j$ consists of
distinct lines.  Define cubics
\[Q_1 = L_{1,2} \cdot L_{3,4} \cdot L_{5,6}, \ Q_2 = L_{1,2} \cdot L_{3,5}
  \cdot L_{4,6},\]
\[
Q_3 = L_{1,5} \cdot L_{2,6} \cdot L_{3,4}, \ Q _4 = L_{1,3} \cdot L_{2,6} \cdot L_{4,5}.\] If they were linearly dependent, then we could find $a,b,c,d \in k$, not all zero, such that the equation
\[a Q_1 + b Q_2 = c Q_3 +d Q_4\] would be satisfied.  But then $L_{1,2}$ divides the lefthand side, so it must also divide the righthand side.  The righthand side is also divisible by $L_{2,6},$ so if it is nonzero, it factors as a product of 3 linear forms.  However, the third form would have to vanish on $p_3, p_4,$ and $p_5$, which contradicts our assumption that the points are in linearly general position. We see that $a = b =c =d=0$ and conclude that the 4 cubics generate $I_{\AAA}.$ 
\end{proof}

For $r > 6$ an elementary dimension count shows that $I_{\mathcal{A}}$
cannot be pl-generated if ${\mathcal{A}}$ consists of $r$ points in linearly
general position.

\begin{prop}\label{prop: points}
Let $\mathcal{A}$ be an arrangement of $r>6$ points in ${\bf P}^2$ in linearly
general position.  Then  $I_{\mathcal{A}}$ is not
pl-generated.
\end{prop}

\begin{proof}
Note that any linear form defines a line in ${\bf P}^2$ that contains at most
two points of $\mathcal{A}.$  Thus, the minimum degree of a product of
homogeneous linear forms that vanishes on $\mathcal{A}$ is $\lceil r/2 \rceil.$
    Hence, we are done if we can show that there must be a form of degree less
than $\lceil r/2 \rceil$ in $I_{\mathcal{A}}.$  This follows from the fact that
   \[r < \binom{\lceil r/2 \rceil - 1+2}{2} = \binom{\lceil r/2 \rceil +1}{2} \]
if $r > 6.$
\end{proof}

\begin{rmk}\label{rmk: points}
An analogous argument shows that for $r \gg 0,$
    the ideal of a linearly general arrangement of $r$ points in ${\bf P}^q$ cannot be
pl-generated.
\end{rmk}

The following proposition gives an example of hypotheses on the geometry of the
   points that imply that $I_{\mathcal{A}}$ is pl-generated.

\begin{prop}\label{prop: pointsonconic}
If ${\mathcal{A}}$ is a set of $r$ points
contained in a union of 2
    lines then $I_{{\mathcal{A}}}$ is pl-generated.
\end{prop}

\begin{proof}
As mentioned at the beginning of the section, any set of 4 points in ${\bf P}^2$ is pl-generated.  So we may assume that $r > 4.$ If there exist lines $L_1$ and $L_2$ containing ${\mathcal{A}},$ and $L_1 \cap L_2$ is a point of ${\mathcal{A}},$ then we can find products of linear forms generating $I_{\mathcal{A}}$ via a construction of Geramita, Gregory and Roberts \cite{geramita-gregory-roberts}, discussed in Chapter 3 of \cite{eisenbud}.

Otherwise, we may assume that the conic $L_1 \cup L_2$ is unique and that $L_1$ contains points 
$p_1, \ldots, p_{r_1} \in {\mathcal{A}}$ and $L_2$ contains $q_1, \ldots, q_{r_2} \in {\mathcal{A}}$ with
$r_1 \ge r_2.$  For $i = 1, \ldots, r_2,$ let $h_i$ be a linear form defining the line joining $p_i$ to $q_i.$  
For $i = r_{2+1}, \ldots, r_1,$ pick any line through $p_i$ not equal to $L_1$ and let $h_i$ be its defining equation.

Since $L_1 \cap L_2 \not\in \AAA,$ the points are a complete
  intersection if $r_1 = r_2,$ and $I_{{\mathcal{A}}}= (L_1 \cdot L_2, h_1 \cdots h_{r_1}).$  If $r_1 > r_2,$ then
 using Corollary \ref{cor: numgens} and Theorem \ref{thm: invariants}
we see that $I_{{\mathcal{A}}}$ must be minimally generated by the conic $L_1 \cdot L_2$
plus generators of degrees $r_2+1$ and $r_1.$  (See also Exercises 3.4--3.7 in \cite{eisenbud}.)
Therefore, 
$$I_{\AAA} = (L_1 \cdot L_2,\  L_1 \cdot h_1 \cdots h_{r_2}, h_1 
\cdots h_{r_1}).$$
\end{proof}

Arrangements of points in ${\bf P}^2$ whose ideals are
pl-generated also appear in \S 2 of \cite{geramita-maroscia}, 
and results in \cite{geramita-gregory-roberts} show that for every Hilbert function of points in ${\bf P}^2$ there exists a finite set of points whose ideal is pl-generated having that Hilbert function, as long as $k$ is infinite. In fact, as Theorem 3.13 in \cite{eisenbud} shows, one can specify the $e_i$ and $f_j$ appearing in Theorem \ref{thm: invariants}.
\bigskip

\subsection{Lines in ${\bf P}^3$}\label{sec: lines}
In this section we will explore when the ideal of an arrangement of lines in
${\bf P}^3$ is pl-generated.

It is easy to construct
line arrangements whose defining ideals are not pl-generated.
Recall that ${\bf P}^{1} \times {\bf P}^{1}$ can be embedded into
${\bf P}^{3}$ as an
    irreducible quadric surface $Z.$ Then, over an infinite field, $Z$ has two
infinite rulings
of disjoint lines $\{X_{\alpha}\}$
and
$\{Y_{\beta}\}.$  Let ${\mathcal{A}}$ be a line arrangement consisting of
$r$ distinct lines from among the set $\{X_{\alpha}\}.$
    Since each pair
of distinct lines in ${\mathcal{A}}$ is disjoint, no two lines are
contained in a hyperplane.
Therefore, the minimum degree of a product of homogeneous linear forms that
vanishes on ${\mathcal{A}}$ is $r.$  Thus,
if $r > 2$, it is clear that $I_{\mathcal{A}}$ cannot
be pl-generated.

More generally, we have the following result.

\begin{prop}  \label{prop: lines}If ${\mathcal{A}}$ is any collection of $r>2$
disjoint lines in ${\bf P}^{3}$,
then $I_{\mathcal{A}}$ cannot be pl-generated.
\end{prop}

\begin{proof}
Since any three skew lines in ${\bf P}^{3}$ lie on an irreducible quadric
surface  (see \cite{harris}, Ex. 2.12.),
it follows  that there is a form of degree 
$2\lfloor r/3 \rfloor +
a$ in $I_{{\mathcal{A}}},$ where $a \equiv r \pmod{ 3}.$  But no pair of the lines 
is contained in a hyperplane, so the
minimum degree of a product of linear forms vanishing on the $r$ lines is $r.$
\end{proof}

\begin{rmk}
Since one does not expect lines in ${\bf P}^3$ to meet, Proposition 
\ref{prop: lines} implies that in practice, lines in ${\bf P}^3$ picked
at random will not have a pl-generated ideal.
\end{rmk}

\begin{rmk}\label{rmk: lines}
Proposition \ref{prop: lines} is a generalization of the example 
given in \cite{li-li2} of three subspaces
  whose ideal is not pl-generated; the 
subspaces given there are in
fact three skew lines in ${\bf P}^3.$  One can generalize the 
statement further to show
that the ideal of $r$ skew $(k-1)$-planes in
  ${\bf P}^{2k-1}$ is not pl-generated when $r \gg 0.$  Three
pairwise disjoint $(k-1)$-planes lie on a variety defined
by quadrics which is projectively equivalent to a Segre variety.
(See again \cite{harris}, Ex. 2.12.)  However, no two of the
$(k-1)$-planes can lie in a hyperplane.
Using products of quadrics, when $r$ is large we can find a form of degree $<r$
  that vanishes
on the arrangement.
\end{rmk}

We close this section with an example showing that there are line arrangements
in ${\bf P}^{3}$ that cannot be scheme-theoretically defined by any 
pl-generated
ideal.

Recall that the \emph{saturation} of a homogeneous ideal $I$ in
$k[x_1, \ldots x_n]$ is
   defined to be \[ \{ f \in k[x_1, \ldots, x_n] \mid f \cdot (x_1, 
\ldots, x_n)^d
  \subseteq I \,
\mathrm{for} \, d \gg 0\}.\]  The saturation of an ideal $I$ is the largest
ideal defining the projective subscheme defined by $I,$ and ideals with
distinct saturations define distinct schemes.

Let $\AAA$ be an arrangement of lines in ${\bf P}^3.$  Construct an embedding
of $\AAA$ into a hyperplane arrangement $\tilde{\HHH}$ as in 
Construction \ref{con: embed}.  Theorem \ref{prop: full ideal} states 
that
$F_{\AAA, \tilde{\HHH}}$ is the
ideal generated by all products of linear forms in $I_{\mathcal{A}}.$  The
following example shows that it may be the case that the three ideals $I_1
\cdots I_r,$ $F_{\AAA, \tilde{\HHH}},$ and $I_{{\mathcal{A}}}$ define 
three different
schemes.  (Thanks to D. Eisenbud and R. Lazarsfeld for suggesting to
investigate
   cones of subspaces.)

\begin{prop}
Let $k[w,x,y,z]$ be the coordinate ring of ${\bf P}^3.$  Let $X$ be a set of
five points in ${\bf P}^2$ in linearly general
position, and let $I_1, \ldots, I_5 \subseteq k[x,y,z]$ be their defining
  ideals.  Let $\tilde{I_i} = I_i \cdot k[w,x,y,z],$  and define $I_{\AAA}
=  \tilde{I_1} \cap \cdots \cap \tilde{I_5},$ so that $\AAA$ is a cone
over $X.$  Then 
$\tilde{I_1}\cdots
\tilde{I_5},$  $F_{\AAA, \tilde{\HHH}},$ and
$I_{{\mathcal{A}}}$
are saturated and are all different.
\end{prop}\label{prop: cones}

\begin{proof}
The ideals
$\tilde{I_1} \cdots \tilde{I_5},$ $F_{\AAA, \tilde{\HHH}},$ and 
$I_{{\mathcal{A}}}$
are all generated
   by
polynomials in $x,y,z$
and are hence
saturated as ideals in $k[w,x,y,z].$

Since each pair of lines lies in a hyperplane, $F_{\AAA, 
\tilde{\HHH}}$ contains
elements of degree three. This shows that it cannot be equal to 
$\tilde{I_1}\cdots
\tilde{I_5}$, which contains only elements of degree $\ge 5.$

Additionally, $F_{\AAA, \tilde{\HHH}}$ cannot be equal to 
$I_{{\mathcal{A}}}$ because
$I_{{\mathcal{A}}}$
    contains the equation of the cone over the unique conic determined
by the points in $X$, but
$F_{\AAA, \tilde{\HHH}}$ contains only forms of degree $\ge 3.$
\end{proof}

\smallskip

\vglue .3cm
\bibliographystyle{amsalpha}

\end{document}